\documentclass[12pt]{article}
\usepackage[cp1250]{inputenc}
\usepackage{amsmath,amstext,amssymb,amsopn,amsthm}

\theoremstyle{plain}
\newtheorem{thm}{Theorem}[section]

\newtheorem{lem}[thm]{Lemma}
\newtheorem{prop}[thm]{Proposition}
\newtheorem{conj}{Conjecture}[section]
\newtheorem{cor}[thm]{Corollary}

\theoremstyle{definition}

\theoremstyle{remark}
\newtheorem*{rem*}{Remark}

\newcommand{\R}{\mathbb{R}}
\newcommand{\N}{\mathbb{N}}

\newcommand{\C}{\mathbb{C}}
\renewcommand{\H}{\mathbb{H}}
\newcommand{\D}{\mathbb{D}}

\renewcommand{\leq}{\leqslant}
\renewcommand{\geq}{\geqslant}

\newcommand{\pref}[1]{(\ref{#1})}

\def\o{\over}
\def\({\left(}
\def\){\right)}
\def\[{\left[}
\def\]{\right]}
\def\<{\langle}
\def\>{\rangle}

\title {Poisson kernel and Green function of the ball in real hyperbolic spaces
 \footnotetext{2000 MS Classification:
                   Primary 60J65; Secondary 60J60.
   {\it Key words and phrases}: hyperbolic spaces, Brownian
   motion, Poisson kernel, Green function.
   Research partially supported by KBN
        grant 1 P03A 020 28 and RTN Harmonic Analysis and Related Problems
        contract HPRN-CT-2001-00273-HARP}}
 \author{ T. Byczkowski and  J. Ma\l{}ecki, \\ Institute of Mathematics, Wroc\l{}aw
    University of Technology,\\ ul. Wybrze\.ze Wyspia\'nskiego 27, 50-370 Wroc\l{}aw, Poland}
\date{}
\begin{document}
\maketitle
\begin{abstract}
Let $(X_t)_{t\geq0}$ be the $n$-dimensional hyperbolic Brownian motion, that is the
diffusion on the real hyperbolic space $\D^n$ having the
Laplace-Beltrami operator as its generator. The aim of the paper is to derive the formulas for the Gegenbauer transform of the
Poisson kernel and the Green function of the ball for the process $(X_t)_{t\geq0}$. Under some additional hypotheses we give the formulas
for the Poisson kernel itself. In particular, we provide formulas in
$\D^4$ and $\D^6$ spaces for the Poisson kernel and the Green function
as well.
\end{abstract}
\newpage

 \section{Introduction}

  Investigation of the  hyperbolic Brownian motion is an important
  and intensely developed topic in recent years (cf. \cite{revista}, \cite{bougerol}).
   On the other hand, it is well known that the Poisson kernel for a region is a
   fundamental tool in  harmonic analysis or probabilistic
 potential theory.  In the classical situation of the Laplacian in $\R^n$,
 the exact formula for the kernel leads to many
 important results concerning behaviour of harmonic functions.
 Moreover, probabilistic potential theory uses Poisson kernel
 techniques to find solutions to
 the Schr\"odinger equation (\cite{ChZ}).
 Availability of the exact formula for the kernel is often of crucial
 importance for the argument.

 In the case of half-spaces in the model $\H^n$ of real hyperbolic spaces
 (or, equivalently, the region bounded by a horocycle, in $\D^n$)
 the form of the Fourier transform of the corresponding Poisson kernel
 is known for some time (see \cite{Du} and \cite{BTF}). However, in many
 applications the resulting Fourier-Hankel inversion formula is of little
 use. A satisfactory integral representation of the Poisson kernel in this
 case, along with the resulting analysis of the asymptotic behaviour,
  was given in \cite{BGS}.

 The aim of this paper is to provide the Gegenbauer transform for the
 Poisson kernel and Green function of a ball in the real hyperbolic space $\D^n$.
 Next, we determine the formula for the Poisson kernel itself; first - as a
 series representation, and, in some cases - as an explicit integral formula.
 Results presented here are not that complete as in \cite{BGS}; they depend
 on the properties of a hypergeometric function $F_k$, which appears
 quite naturally in the Gegenbauer transform of the Poisson kernel.

 Although we are motivated here by the paper
 \cite{W}, where the Gegenbauer transform of the
 joint distribution of hitting time and hitting distributions for the ball in
 the case of classical Brownian motion in $\R^n$ was found, there were very
 substantial difficulties to adapt the above approach to the present situation.

 We also provide explicit formulas for Poisson kernel and Green function
 in $\D^4$ and $\D^6$.

  The paper is organized as follows. In Section 2, after
  some preliminaries, we apply stochastic calculus to write a "polar" decomposition
  of the hyperbolic Brownian motion on $\D^n$.
  It is the starting point to obtain, in Section 3, the basic formula for
  Gegenbauer coefficients of the cosine between the axis determined by the
  process and the starting point $x\neq 0$ (see Theorem \ref{thm01}).
  We again use here the
  stochastic calculus; more specifically Feynmann-Kac technique and the relation
  with the appropriate Schr\"odinger equation. The series representation theorem
  of the Poisson kernel is the main result of this section.

  In Section 4, in Theorem \ref{thm03}, we provide Gegenbauer coefficients for
  the Green function of the ball, restricted to a suitable sphere. In the next theorem
  we write series representation theorem for the Green function of the ball. Remark that
  although the Poisson kernel determines uniquely the Green function (by sweeping out
  formula), the substantial computational complexity forces us to find an alternative
  approach.

  In Section 5 we try to describe an explicit integral representation of Poisson kernel.
  It turns out that all depends on the properties of the hypergeometric function $F_k$.
  Conditions sufficient to obtain the desired representation are collected in
   Conjectures \ref{conjecture01} and \ref{conjecture02}.
  Checking the conditions imposed on the specific hypergeometric functions $F_k$
  we exhibit the exact form of this representation for $\D^4$ and $\D^6$. In general,
  the validity of this conjecture depends on the location of roots of $F_k$, with
  respect to the (complex) variable $k$, which seems to be a difficult problem.
   We also provide explicit formulas for
  the Green function of the ball for $\D^4$ and $\D^6$.

 \section{Preliminaries}

We begin with some basic informations about hypergeometric
 functions and Gegenbauer polynomials, needed in the sequel. This part of the material is
 standard and can be found, e.g. in \cite{E}. After that we identify Brownian
 motion in real hyperbolic spaces, in terms of Stochastic Differential Equations
 (SDE).
 We then discuss briefly properties of the heat kernel on hyperbolic spaces, following
 approach presented in the monograph of E. B. Davies \cite{D}. In the end we obtain a kind
 of "polar" decomposition of the hyperbolic Brownian motion, in terms of SDE.

We denote here by $(x,y)$ the standard inner product of $x,y \in \R^n$
and by $|x|$ the Euclidean length of a vector $x$. The
sphere with center at $0$ and the radius $r$ is written as $S_r = \{x\in R^n:
|x|=r\}$.  The $(n-1)$-dimensional spherical measure on  $S_r$ will by
denoted be $\sigma_r$.
 Put
\begin{equation*}
    \omega_{n-1} = \frac{2\pi^{n\o2}}{\Gamma({n\o2})}, \quad
    n=1,2,\ldots
\end{equation*}
It is the total mass of the associated $(n-1)$-dimensional spherical measure of
the unit sphere $S_1$. Note that for $n=1$, $S_1$ is a two-point
set and its $0$-dimensional measure is equal to the counting measure by an accepted
convention. For the rest of the paper we assume that $n>2$.

We will denote by $F(\alpha,\beta;\gamma;z)$ the hypergeometric
function of variable $z$ with parameters $\alpha,\beta,\gamma$.
For $|z|<1$ and $\alpha,\beta,\gamma \in \C$, $\gamma \neq
0,-1,-2,\ldots$ the function $F$ is defined by the hypergeometric
series
\begin{equation*}
F(\alpha,\beta;\gamma;z) = \sum_{i=0}^\infty
\frac{(\alpha)_i(\beta)_i}{(\gamma)_i i!}\,z^i \/,
\end{equation*}
where $(\alpha)_i = \Gamma(\alpha+i)/\Gamma(\alpha)$ is the
Pochhammer symbol. We shall supplement the definition in the case
$\alpha=-l$ and $\gamma=-m$ where $l=0,1,2,\,\ldots\,$, and
$m=l,l+1,l+2,\,\ldots\,$. Then it is customary to define
\begin{equation*}
F(-l,\beta;-m;z) = \sum_{i=0}^l \frac{(-l)_i(\beta)_i}{(-m)_i
i!}\,z^i \/.
\end{equation*}
To simplify our notation we put $\rho=\frac{n-2}{2}$ and define
\begin{equation*}
F_k(z)=F(k,-\rho;k+{n\o2};z) \/,
\end{equation*}
and for $k>0$ or $n/2 \notin \N$
\begin{equation*}
G_k(z)=F(-\rho,2-k-n;2-k-{n\o2};z) \/.
\end{equation*}
When $k=0$ and $n$ is an even number greater than $2$ we put
\begin{equation*}
G_0(z)= \rho\,\sum_{i=0,\,\,i\neq\rho}^{n-2}\binom{n-2}{i}
\frac{(-1)^{i+1}}{i-\rho}z^i
+\rho\,\binom{n-2}{\rho}(-1)^{n\o2}\,z^\rho\log z \/.
\end{equation*}
The general solution of the hypergeometric equation
\begin{equation*}
z(1-z)u''+(\gamma-(\alpha+\beta+1)z)u'-\alpha\beta u = 0
\end{equation*}
for $\alpha=k, \beta=-\rho$ and $\gamma=k+{n\o2}$ is given by
\begin{equation}
\label{general_solution}
c_1\cdot
F_k(z)+\frac{c_2}{z^{k+\rho}}\cdot G_k(z) \/,
\end{equation}
where $c_1,c_2$ are constants and $k=0,1,2,\,\ldots\,$. Observe that
$F_k(z)$ is bounded and $z^{-k-\rho}G_k(z)$ (the second solution of
(\ref{general_solution})) is unbounded on the interval $(0,a]$ for
every $a \in (0,1)$. It follows from the fact that the functions
$F_k(z),\,G_k(z)$ are continuous on $[0,a]$.

We also have
\begin{prop}
\begin{equation}
\label{hypgeo01}
(k+\rho)F_k(z)G_k(z)+zF_k'(z)G_k(z)-zF_k(z)G_k'(z)=(k+\rho)(1-z)^{n-2}.
\end{equation}
\end{prop}
\begin{proof}
Let denote by $u(z)$ the function on the left-hand side of
(\ref{hypgeo01}). Using the hypergeometric equations which are
satisfied by $F_k$ and $G_k$ we find that $(1-z)u'(z)=(2-n)u(z)$.
We can also show that $u(0)=k+\rho$. So the desired equality
follows. The details are left to the reader.
\end{proof}

Gegenbauer's polynomial $C_k^{(v)}(z)$ for integer value of $k$
and $v>0$ is defined to be the coefficient of $h^k$ in the
Maclaurin expansion of $(1-2zh+h^2)^{-v}$, considered as a function of $h$.
So we have
\begin{equation*}
    (1-2zh+h^2)^{-v} = \sum_{k=0}^{\infty}C_k^{(v)}(z)h^k, \quad
    |z|\leq 1, |h|<1.
\end{equation*}
Observe that $C_0^{(v)}(z)=1$, for all $v>0$.
For $v=0$ it is customary to take $C_0^{(0)}(z)\equiv 1$,
$C_k^{(0)}(z)=\lim_{v \to
0}\frac{C_k^{(v)}(z)}{v}=\frac{2T_k(z)}{k}$, where $T_k$ is $k$th
Chebyshev polynomial defined by $T_k(\cos\phi)=\cos(n\phi)$. One of
generating functions for $T_k$ is given by
\begin{equation*}
\log(1-2zh+h^2)^{-1} = 2\sum_{k=1}^{\infty}k^{-1}T_k(z)h^k.
\end{equation*}
We also have the following trigonometric expansion of
$C_k^{(v)}(\cos\phi)$:
\begin{equation*}
\Gamma(v)^2C_k^{(v)}(\cos\phi) =
\sum_{l=0}^k\frac{\Gamma(l+v)\Gamma(k-l+v)}{l!(k-l)!}\cdot
e^{-i(k-2l)\phi},
\end{equation*}
which gives
\begin{equation}
\label{gegen01}
    |C_k^{(v)}(\cos\phi)|\leq C_k^{(v)}(1), \quad \phi \in [0,\pi].
\end{equation}
Note that $C_k^{(v)}(1)=\Gamma(k+2v)/(k!\,\Gamma(2v))$. We recall
the orthogonal relations of Gegenbauer polynomials
\begin{eqnarray*}
\int_{-1}^{1}C_k^{(v)}(x)C_l^{(v)}(x)(1-x^2)^{v-{1\o2}}dx =
\delta_{kl}\frac{2^{1-2v}\pi\Gamma(k+2v)}{k!(v+k)\Gamma(v)^2}.
\end{eqnarray*}
In the special case for $v=\rho$, using the well known relation
for gamma function
\begin{equation*}
2^{2z-1}\Gamma(z)\Gamma(z+{1\o2})=\sqrt{\pi}\,\Gamma(2z),
\end{equation*}
we obtain
\begin{equation}
\label{gegen02}
\int_{-1}^{1}C_k^{(\rho)}(x)C_l^{(\rho)}(x)(1-x^2)^{{n-3}\o2}dx =
\delta_{kl}\frac{\rho}{k+\rho}\cdot C_k^{(\rho)}(1)\cdot
\frac{\omega_{n-1}}{\omega_{n-2}}.
\end{equation}

Note also that the polynomial $C_k^{(v)}(z)$ is a solution of the
Gegenbauer differential equation
\begin{equation}
\label{gegen03}
    (z^2-1)\omega''+(2v+1)z\omega'-k(k+2v)\omega = 0.
\end{equation}
Now, we introduce some basic information about measures
(functions) on spheres and their Gegenbauer transforms. For more
details see \cite{W} and \cite{E}.

We say that a finite Borel measure $\mu(\cdot)$ on sphere $S_r$ is
$axially\,symmetric$ (AS) with axis $x\in\R^n$ if $\mu(UA)=\mu(A)$
for each Borel set $A$ and each orthogonal transformation $U$ such
that $Ux=x$. The definition of AS functions is similar. For AS
measures we define its Gegenbauer coefficient $\widehat{\mu}_k$ by
\begin{equation*}
C_k^{(\rho)}(1)\widehat{\mu}_k =
\int_{S_r}C_k^{(\rho)}(\cos\theta)\mu(dy) \/,
\end{equation*}
where $\theta = \angle x0y$ if $x \neq 0$ and $\theta = \angle
u0y$, for arbitrary but fixed nonzero vector $u$, in case $x=0$.
We will use the following property of Gegenbauer transform:
\begin{thm}
\label{thm0} The AS measure $\mu$ is uniquely determined by its
transform $\{\widehat{\mu}_k\}_{k=0}^{\infty}$.
\end{thm}
Clearly, the same is true for AS functions, i.e. any AS function is
uniquely determined by its Gegenbauer transform.

Consider the ball model of the $n$-dimensional real hyperbolic
space
\begin{equation*}
\D^n = \{x\in \R^{n} : |x|<1 \}, \quad n>2 .
\end{equation*}
The Riemannian metric and the distance formula are given by
\begin{equation}
    \label{rmetric}
  ds^2 = {|dx|^2 \o (1-|x|^2)^2} \, ,
\end{equation}
\begin{equation*}
    \cosh(2d(x,y)) = 1+\frac{2|x-y|^2}{(1-|x|^2)(1-|y|^2)}\,.
\end{equation*}
The canonical (hyperbolic) volume element is given by
\begin{equation*}
dV_n={ dx \o  (1-|x|^2)^{n}} \,.
\end{equation*}

Consider the following system of SDE:
\begin{equation*}
    \label{sde01}
    \frac{dX_k(t)}{1-|X(t)|^2} = dB_k(t)+2(n-2)X_k(t)dt, \quad
    k=1,\ldots,n \/.
\end{equation*}
For $f \in \mathcal{C}^2$ by It\^o Formula we obtain
\begin{eqnarray*}
    f(X_t)-f(X_0) &=&
    \sum_{i=1}^n \int_0^t \dfrac{\partial f}{\partial x_i}dX_i(s)+{1 \o
    2}\sum_{i=1}^n \int_0^t \dfrac{\partial^2f}{\partial
    {x_i}^2}d\langle X_i\rangle (s)\\
    &=& \sum_{i=1}^n\int_0^t \dfrac{\partial f}{\partial
    x_i}(1-|X(s)|^2)dB_i(s)+2(n-2)\sum_{i=1}^n \int_0^t\dfrac{\partial f}{\partial
    x_i}(1-|X(s)|^2)X_i(s)ds\\
    &+&\sum_{i=1}^n \int_0^t\dfrac{\partial ^2f}{\partial
    x_i^2}(1-|X(s)|^2)^2ds  \/.
\end{eqnarray*}
Here $B=(B_1,\ldots,B_n)$ denotes the standard Brownian motion
with scaling such that $EB_k^2(t) = 2t$, $k=1,\ldots,n$. Thus, the
generator of the process $X(t)$ determined by the above system of SDE
is given  by
\begin{equation*}
    \Delta_B = (1-|x|^2)^2\sum_{i=1}^n \dfrac{\partial^2}{\partial
    x_i^2}+2(n-2)(1-|x|^2)\sum_{i=1}^n x_i \dfrac{\partial}{\partial x_i} \/,
\end{equation*}
and is the canonical Laplace - Beltrami operator associated with the
Riemannian metric (\ref{rmetric}).
Since the Laplace - Beltrami operator commutes with isometries acting
on $\D^n$, the heat kernel $k_n$, i.e. the transition density of the hyperbolic Brownian
motion is a function of the (hyperbolic) distance $d(x,y)$. Fix $a \in \D^n$ and
denote $\rho(x)=d(a,x)$.
We have the following explicit form of the heat kernel $k_2(t,\rho)$ on the hyperbolic
disc $\D^2$ (see \cite{D}):
\begin{equation*}
k_2(t,\rho)=2^{3\o2} (4\pi t)^{-3/2} e^{-t} \int_{\rho}^\infty {s
e^{-s^2/4t} \,ds \o (\cosh (2s) - \cosh (2\rho))^{1/2}} \, ,
\end{equation*}
while on $\D^3$ we obtain
\begin{equation*}
k_3(t,\rho)= 2(4\pi t)^{-3/2} e^{-4t - \rho^2/4t} {\rho \o \sinh
(2\rho)} \,.
\end{equation*}
In higher dimensions we have
the following recursion formula:
\begin{equation} \label{heatker}
k_n(t,\rho)= \sqrt{2} e^{(2n-1)t} \int_{\rho}^\infty
{k_{n+1}(t,\lambda) \sinh{(2\lambda)} d\, \lambda \o
(\cosh{(2\lambda)} - \cosh{(2\rho}))^{1/2}} \,.
\end{equation}
We now show that the hyperbolic Brownian motion on $\D^2$ is {\it transitive}
for $n\geq 2$ (see, e. g. \cite{Ch}), a fact which is widely known; we include it
for the reader's convenience. For this purpose it is enough to show
that the potential $U_n(z)<\infty$, for almost all $0<z\in \R$.
A direct computation on $\D^2$ yields
\begin{equation*}
U_2(\rho)=\int_0^\infty k_2(t,\rho) \, dt = \sqrt{2} (4\pi)^{-1} \ln
\coth (\rho) \,.
\end{equation*}

Note that $k_n$ is the density function with respect to the canonical hyperbolic
volume element $dV_n$.

The hyperbolic Brownian motion process on $\D^n$ with the measure $dV_n$
as the reference measure fits into context of the so-called "dual processes"
and their potential theory (see \cite{BG}, ch. VI).
From this theory it follows, in particular, that single points in $\D^n$ are polar,
whenever we show that the potential kernel $U_n(\rho)<\infty$, for $\rho >0$.
Indeed, by Proposition 3.5, Ch. II, in \cite{BG} we obtain that $\{a\}=\{U_n(\rho)=\infty\}$
is polar if and only if it is null (of potential 0). This last condition is obviously
satisfied whenever there exists almost everywhere finite potential kernel $U_n$.
Thus, the Brownian motion on $\D^2$ does not hit single points, a.s. For higher dimensions
we use recursion formula \pref{heatker} to obtain:

\begin{eqnarray*}
\int_{\rho}^\infty {\sinh (2\lambda)\,U_{n+1}(\lambda)  \o (\cosh
(2\lambda) - \cosh (2\rho))^{1/2}} \, d\lambda &=&
\int_{\rho}^\infty {\sinh (2\lambda) \o (\cosh (2\lambda) - \cosh
(2\rho))^{1/2}}
\int_0^\infty k_{n+1}(t,\lambda)\,dt\,d\lambda \\
=2^{-1/2}\int_0^\infty e^{-(2n-1)t} k_{n}(t,\lambda)\,dt &\leq&
 2^{-1/2}\int_0^\infty k_{n}(t,\lambda)\,dt  =2^{-1/2} U_n(\lambda) \,.
\end{eqnarray*}
By induction, we obtain for $n\geq 2$
\begin{equation}\label{transience}
                 U_n(\rho)<\infty \quad \text{a.e.}
\end{equation}
For the hyperbolic Brownian motion $X_t$ let $x=X(0)$ be the
starting point. We define the following two processes :
  \begin{equation} \label{decomp}
 R_t  =
\sum_{i=1}^n X_i^2(t) = |X_t|^2\/, \qquad \Phi_t=\cos\angle x0X_t =
{(x,X(t)) \o |x| |X(t)|} \/.
  \end{equation}
We always assume that $x \neq 0$. Then also $X(t) \neq 0$, a.s. (see \pref{transience}).
Thus, $\Phi_t$ is well defined almost surely.

\begin{prop}
The process $(R(t),\Phi(t))$ satisfies the following system
of stochastic differential equations:
\begin{equation}
    \label{sde02}
  \left\{
    \begin{array}{ccc}
      dR(t) & = & 2(1-R(t))\left(\sqrt{R(t)}dW_1(t)+((n-4)R(t)+n)dt\right) \\
      d\Phi(t) & = & (1-R(t))\left(\frac{1-\Phi^2(t)}{R(t)}\right)^{1 \o
    2}dW_2(t)-(n-1)\frac{(1-R(t))^2}{R(t)}\Phi(t)dt,
    \end{array}\right.
\end{equation}
where $W_1(t),W_2(t)$ are independent Brownian motions
on $\R$ with variation $2t$.
\end{prop}
\begin{proof}
    We define
    \begin{eqnarray*}
    W_1(t) & = & \int_0^t \frac{\sum_{i=1}^n
    X_i(s)dB_i(s)}{\sqrt{R(s)}},\\
    W_2(t) & = & \int_0^t
    \textbf{1}_{(\Phi(s)<1)} \left(\frac{R(s)}{1-\Phi^2(s)}\right)^{1 \o 2}
    \sum_{i=1}^n \left( \frac{x_i}{|x||X(s)|}-
    \frac{X_i(s)(x,X(s))}{|x||X(s)|^3}\right)dB_i(s)\\
    &+& \int_0^t\textbf{1}_{(\Phi(s)=1)}d\widetilde{B}(s),
\end{eqnarray*}
 where $\widetilde{B}(t)$ is classical
 Brownian motion on $\R$ (with variation $2t$) such that
 $B_1(t),\ldots,B_2(t),\widetilde{B}(t)$ are
 independent. It is clear from the definitions and property of the It\^o integral
 that $W_1(t)$, $W_2(t)$ are local martingales. We also have
 \begin{eqnarray*}
    d\<W_1,W_1\>(t) &=& \frac{\sum_{i=1}^n X_i^2(t)}{R(t)}2dt = 2dt \/, \\
    d\<W_2,W_2\>(t)
    &=&
    \textbf{1}_{(\Phi(s)<1)}\frac{R(t)}{1-\Phi^2(t)}
    \sum_{i=1}^n\left(\frac{x_i^2}{|x|^2|X(t)|^2}-
    \frac{2x_iX_i(t)\(x,X(t)\)}{|x|^2|X(t)|^4}\right.\\
    &+&
     \left. \frac{X_i^2(t)\(x,X(t)\)^2}{|x|^2|X(t)|^6}\right)2dt
    + \textbf{1}_{(\Phi(t)=1)}2dt = 2dt \/, \\
    d\<W_1,W_2\>(t)
    &=&
    \textbf{1}_{(\Phi(t)<1)}\left(1-\Phi^2(t)\right)^{-{1 \o
    2}}\sum_{i=1}^{n}\left(\frac{x_iX_i(t)}{|x||X(t)|}-\frac{X_i^2(t)\(x,X(t)\)}{|x|
    |X(t)|^3}\right)2dt = 0 \/.
\end{eqnarray*}
So $W_1(t),W_2(t)$ are independent Brownian motions on
$\R$. Using (\ref{sde01}) and It\^o Formula we get
\begin{eqnarray*}
    d R(t) &=& \sum_{i=1}^n 2X_i(t)dX_i(t)+\sum_{i=1}^n
    d\<X_i\>(t) \\
    &=& 2(1-|X(t)|^2)\left(\sum_{i=1}^n
    X_i(t)dB_i(t)+2(n-2)\sum_{i=1}^n
    X_i^2(t)dt+ n(1-|X(t)|^2)dt\right) \\
    &=& 2(1-R(t))\left(\sqrt{R(t)}\,dW_1(t)+((n-4)R(t)+n)\,dt\right).
\end{eqnarray*}

For the function $g(y)=\frac{(x,y)}{|x||y|}$ we have

\begin{equation*}
\frac{\partial g}{\partial y_i} =
\frac{x_i}{|x||y|}-\frac{(x,y)y_i}{|x||y|^3}, \quad \frac{\partial
^2g}{\partial y_i^2} =
-2\frac{x_iy_i}{|x||y|^3}-\frac{(x,y)}{|x||y|^3}+3\frac{(x,y)y_i^2}{|x||y|^5}
\end{equation*}

and

\begin{equation*}
\sum_{i=1}^n \frac{\partial ^2g}{\partial y_i^2}
=-(n-1)\frac{(x,y)}{|x||y|^3}.
\end{equation*}

 So It\^o Formula gives
\begin{equation*}
    d\Phi(t) = (1-R(t))\left(\frac{1-\Phi^2(t)}{R(t)}\right)^{1 \o
    2}dW_2(t)-(n-1)\frac{(1-R(t))^2}{R(t)}\Phi(t)dt.
\end{equation*}
\end{proof}
The following lemma will be useful in the next sections. The proof is standard
and is omitted.

\begin{lem}
\label{lemma01}
 Let $W=(W_1,W_2)$ be the $2$-dimensional Brownian motion.
 Suppose that $\Psi(t)$ is a real measurable process adapted to
$\mathcal{F}(W(t))$ such that for every $t>0$, $E^x[\int_0^t
\Psi^2(s)ds]<\infty$, $x\in \R^2$, and $Y(t)$ be square integrable process
adapted to $\mathcal{F}(W_1(t))$. Then we have
\begin{equation*}
E^x[Y(t) \int_0^t \Psi(s)\, dW_2(s)] = 0 \/.
\end{equation*}
\end{lem}
%
%

\section{Poisson kernel of the ball}

Consider a ball $D=\{x \in \D^n: |x|<r\}$ for some fixed $0<r<1$.
Recall, using the distance formula, that $D$ is a hyperbolic ball
with center at $0$ and the radius $\frac{1}{2}\log \frac{1+r}{1-r}$.
Define
\begin{equation*}
\tau_D = \inf \{ t \geq 0: X(t)\notin D \} = \inf \{t \geq 0 :
|X(t)| \geq r \}.
\end{equation*}
By $P_r(x,y)$, $x\in D$, $y\in \partial D=S_r$ we denote the
Poisson kernel of $D$, i.e. the density of the harmonic measure
defined by $\mu_x(A)=P^x(X(\tau_D) \in A)$ for every Borel subset
$A \subset S_r$ and $x\in D$.

Let $U$ be an orthogonal transformation that leaves the starting
point $x$ fixed. It easy to see that the processes $X_t$ and
$U^{-1}(X_t)$ have the same distribution under $P^x$. Thus, for
each Borel set $A \subset S_r$ we get
\begin{equation*}
P^x(X(\tau_D) \in U(A)) = P^x(U^{-1}(X(\tau_D)) \in A) =
P^x(X(\tau_D) \in A).
\end{equation*}
Consequently, the harmonic measure is AS measure with axis $x$. It
follows that the Poisson kernel $P_r(x,\cdot)$, as the density function
of the harmonic measure, is AS function on
$S_r$ with axis $x$.

Below we identify the Poisson kernel in terms of its Gegenbauer
transform. We give explicit formula for its Gegenbauer
coefficients.
\begin{thm}\label{thm01}
    For $|x|<r$ we have
\begin{equation}
\label{pkgegenbauer}
 \frac{E^xC_k^{(\rho)}(\Phi(\tau_D))}{C_k^{(\rho)}(1)} =
 \left(\frac{|x|}{r}\right)^k\frac{F_k(|x|^2)}{F_k(r^2)},
\end{equation}
with $k = 0,1,2,\ldots$.
\end{thm}

\begin{proof}
We define
\begin{eqnarray}
    \nonumber
    V(t) &=& \exp(\int_0^t q(R(s))ds) ,\\
    \label{Z_process}
    Z(t) &=& \varphi (\Phi(t))V(t) =
    C_k^{(\rho)}(\Phi(t))\exp(\int_0^tq(R(s))ds),
\end{eqnarray}
where $q(x)=k(k+n-2)\frac{(1-x)^2}{x}$ and $\varphi(z)=C_k^{(\rho)}(z)$
is the Gegenbauer polynomial introduced in Preliminaries.
Observe that the process $V$ has bounded variation sample paths so we
obtain
\begin{equation*}
dV(t) = q(R(t))V(t)dt\/.
 \end{equation*}
The It\^o formula yields:
\begin{eqnarray*}
    dZ(t) &=& \varphi^{'}(\Phi_t)V(t)d\Phi(t)+ \varphi(\Phi_t)dV(t)+{1 \o
    2}\, \varphi^{''}(\Phi_t)V(t)d\<\Phi\>(t) \\
    &=&  \varphi^{'}(\Phi_t)V(t)(1-R(t))\left(\frac{1-\Phi^2_t}{R(t)}\right)^{1 \o
    2}dW_2(t)\\
    &&+ V(t)\frac{(1-R(t))^2}{R(t)}\left((1-\Phi^2_t) \varphi''(\Phi_t)-(n-1)
    \Phi_t \varphi'(\Phi_t)
    +k(k+n-2) \varphi(\Phi_t)\right)\/dt\/.
\end{eqnarray*}

Using (\ref{gegen03}) we obtain
\begin{equation*}
    dZ(t) =
    \varphi^{'}(\Phi_t)V(t)(1-R(t))\left(\frac{1-\Phi^2_t}{R(t)}\right)^{1 \o 2}
    dW_2(t) \/.
\end{equation*}
Denote


\begin{equation*}
T_n=\inf\{t>0; R(t)\leq 1/n \}\/.
\end{equation*}
Since $R(t)>0$ we obtain $T_n \to \infty $.

From the definition of $Z_t$ and the previous equality we obtain
\begin{eqnarray*}
&{}&    E^xC_k^{(\rho)}(\Phi(t \land T_n \wedge \tau_D )) =
 E^x\left(Z(t \land T_n  \wedge \tau_D)
    \exp(-\int_o^{t \land T_n \wedge \tau_D} q(R(s))ds)\right)\\
    &=& C_k^{(\rho)}(1)\,E^{|x|^2}\exp(-\int_0^{t \land T_n \wedge\tau_D}
    q(R(s))ds)+E^x\left(V^{-1}(t \land T_n \wedge \tau_D)\int_0^{t\wedge T_n}
      \Psi(s)\,dW_2(s)\right) \/,
\end{eqnarray*}
where
\begin{equation*}
\Psi(t) = \textbf{1}_{\{t\land\tau_D\}}\,
\varphi'(\Phi_t)V(t)(1-R(t)) \left(\frac{1-\Phi^2_t}{R(t)}\right)^{1
\o2}  \/.
\end{equation*}
 Since $\tau_D$ and $T_n$ depend only on $R_t$ (i.e. $W_1$) we can
use Lemma \ref{lemma01} to show that the last expectation is
equal to zero. Moreover, since $t \land T_n \wedge \tau_D$ tends to $\tau_D$
as $t \to \infty$ and $n \to \infty$, and the
 functions $C_k^{(\rho)}(\Phi(t \land T_n \wedge \tau_D)$
and $\exp(-\int_0^{t \land T_n \wedge \tau_D} q(R(s))ds)$ are bounded by a
constant, we obtain (using dominated convergence theorem) that
\begin{equation}
  E^xC_k^{(\rho)}(\Phi(\tau_D)) = C_k^{(\rho)}(1)E^{|x|^2}e_{-q}(\tau_D) \/,
\end{equation}
where $e_{-q}(\tau_D) = exp(-\int_0^{\tau_D} q(R(s))ds)$.

Observe that the function $\phi(y) = E^ye_{-q}(\tau_D)$ is by
definition the gauge function for the Schr\"odinger operator based on the
generator of the process $R_t$ and the potential $(-q)$. According to the general theory
(see \cite{ChZ}) the function $\phi$ is
the solution of the appropriate Schr\"odinger equation. Using (\ref{sde02}) we
obtain the following formula for the generator of $R_t$:
\begin{equation*}
 4(1-x)^2x\dfrac{d^2}{dx^2}+2(1-x)((n-4)x+n)\dfrac{d}{dx}.
\end{equation*}
Hence, $\phi$ satisfies the following equation
\begin{equation}
    \label{equation2}
    4(1-x)^2x\phi''(x)+2(1-x)((n-4)x+n)\phi'(x)-k(k+n-2)\frac{(1-x)^2}{x}\phi(x)=0
\end{equation}
on $(0,r^2]$. Let $\phi(x)=x^{k \o 2}g(x)$. Then $\phi'(x) = {k \o
2}x^{{k-2}\o 2}g(x)+x^{k \o 2}g'(x)$ and $\phi''(x) = {{k(k-2)}\o
4}x^{{k-4}\o 2}g(x)+kx^{{k-2} \o 2}g'(x) + x^{k \o 2}g''(x)$ and
consequently (\ref{equation2}) reads as
\begin{equation}
    \label{equation3}
    x(1-x)g''(x)+(k+{n\o2}-(k-{{n-2}\o2}+1)x)g'(x)+k{{n-2}\o2}g(x)=0.
\end{equation}
This is the hypergeometric equation with parameters $\alpha = k$,
$\beta = -\rho$ and $\gamma = k+{n\o2}$. The general solution of
(\ref{equation3}) is given by (\ref{general_solution}) and we
infer that
\begin{equation*}
    \phi(x) =
    x^{k\o2}\left(c_1 \cdot
    F_k(x)+\frac{c_2}{x^{\rho+k}}G_k(x)\right)\/.
\end{equation*}
By definition $\phi(x)$ is bounded and $\phi(r^2)=1$. Therefore
$c_2=0$ and the other condition gives the normalizing constant
\begin{equation*}
    c_1 = \frac{1}{r^k{F_k(r^2)}}.
\end{equation*}
This completes the proof.
\end{proof}

\begin{thm}[Poisson kernel formula]
\label{thm02}
    For $|x|<r$, $|y|=r$ we have
\begin{equation}
 \label{pkformula}
 P_r(x,y) = \frac{\Gamma({n \o 2})}{2\pi^{n \o 2}r^{n-1}} \sum_{k=0}^{\infty}\frac{k+\rho}
 {\rho}\,
 \left(\frac{|x|}{r}\right)^k\frac{F_k(|x|^2)}{F_k(r^2)}\,\,
 C_k^{(\rho)}(\cos\theta).
\end{equation}

\end{thm}
\begin{proof}
    It is easy to see that the terms of the series (\ref{pkformula}) are continues functions
    of $y$. We will show that the series is uniformly convergent on the sphere $S_r$.
    Indeed, we have for every $|x|<1$
\begin{equation*}
    |F_k(x)| \leq
    \sum_{j=0}^{\infty}\frac{(k)_j}{(k+{n\o2})_j}\frac{|(-\rho)_j|}{j!}|x|^j
    \leq \sum_{j=0}^{\infty}\frac{|(-\rho)_j|}{j!}|x|^j < \infty
\end{equation*}
 and using the above inequality and uniform convergence
 we obtain that
 \begin{equation*}
 \lim_{k \to \infty}F_k(x) =
 \sum_{j=0}^{\infty}\frac{(-\rho)_j}{j!}x^j=(1-x)^{\rho}.
 \end{equation*}
 Moreover, as a consequence of the relation $E^x(e_{-q}(\tau)) =
 \left(\frac{|x|}{r}\right)^k \frac{F_k(|x|^2)}{F_k(r^2)}$, we
 get $F_k(r^2)\neq 0$ for every $k \in \N$.
 Thus there exists a constant $c=c(r,|x|)$ such that
 \begin{equation*}
 \left\vert\frac{F_k(|x|^2)}{F_k(r^2)}\right\vert \leq
 c.
 \end{equation*}
 From this and (\ref{gegen01}) we obtain
 \begin{equation}
    \label{inq01}
    \left\vert\sum_{k=0}^{\infty}\frac{k+\rho}{\rho}\,
    \left(\frac{|x|}{r}\right)^k\frac{F_k(|x|^2)}{F_k(r^2)}\,\,
    C_k^{(\rho)}(\cos\theta)\right\vert < c\sum_{k=0}^\infty \frac{k+\rho}{\rho}
    \left(\frac{|x|}{r}\right)^k C_k^{(\rho)}(1)<\infty  \/.
 \end{equation}
Applying the formula for
$C_k^{(\rho)}(1)$ we obtain
$ C_{k+1}^{(\rho)}(1)/C_k^{(\rho)}(1)$ $=(n+k-1)/(k+1) $. Using
Ratio Criterion for convergence of power series we infer
that the above series is convergent for $|x|<r$.
 Consequently,
the series (\ref{pkformula}) represents continues function of $y$
which is axially symmetric. Using (\ref{gegen02}) and
(\ref{inq01}) we can check that the Gegenbauer coefficient of that
function is the same as the Gegenbauer coefficient of the Poisson
kernel computed in Theorem \ref{thm01}. The desired equality
follows from Theorem \ref{thm0}.

\end{proof}

\section{Green function of the ball}
Let $(X_t^D,P_t^D)$ be the hyperbolic Brownian motion killed
at the boundary $\partial D$. The Green function of $D$ is defined by
\begin{equation*}
G_D(x,y) = \int_0^\infty p_t^D(x,y)dt \quad x,y \in D, \quad x\neq
y.
\end{equation*}
where $p_t^D(x,y)$ is the density function of $P_t^D$.
Similar arguments as in the case of the Poisson kernel show that
\begin{equation*}
P_t^D(U(A)) =P^x(t<\tau_D;X_t \in U(A)) = P^x(t<\tau_D;X_t\in A) = P_t^D(A),
\end{equation*}
for each orthogonal transformation $U$ such that $Ux=x$. Thus, the
Green function $G_D(x,\cdot)$ as a function on the sphere $S_R$,
$0<R<r$, is AS function with axis $x$. Recall that its Gegenbauer
coefficient $\widehat{(G_D)}_k(x,R)$ is defined by
\begin{equation*}
C_k^{(\rho)}(1)\widehat{(G_D)}_k(x,R) =
\frac{1}{\omega_{n-1}R^{n-1}}\int_{S_R}C^{(\rho)}_k(\cos\theta)G_D(x,y)d\sigma_R(y),
\quad x \neq 0.
\end{equation*}
Here $\cos\theta = {(x,y) \o |x| |y|}$.
Observe also that $\widehat{(G_D)}_k(x,R)=\widehat{(G_D)}_k(|x|,R) =
\widehat{(G_D)}_k(R,|x|)$. So from now on  we write
$\widehat{(G_D)}_k(|x|,R)$.
 Moreover, in view of the above-mentioned symmetry, it is
enough to determine   $ \widehat{(G_D)}_k(|x|,R) $
  for $|x|<R$.

\begin{thm}
\label{thm03}
For $|x|<R<r$ we have
\begin{eqnarray*}
\widehat{(G_D)}_k(|x|,R) &=& C_n\,
\dfrac{\rho}{k+\rho}\cdot|x|^kF_k(|x|^2)R^k\left(\dfrac{G_k(R^2)}{R^{2k+2\rho}}
-\dfrac{G_k(r^2)}{r^{2k+2\rho}}\cdot\frac{F_k(R^2)}{F_k(r^2)}\right)\/,
\end{eqnarray*}
where $C_n=\frac{\Gamma({n\o2}-1)}{4\pi^{n\o2}}$.
 \end{thm}

\begin{proof}
Let $(R_t^{\widetilde{D}},P_t^{\widetilde{D}})$ be the process
$R_t$ killed
at the boundary of $\widetilde{D}=\{x\in(0,1) : x<r^2\}$.
Observe that $R_t^{\widetilde{D}}=|X_t^D|^2$ and
$\tau_{\widetilde{D}}=\tau_D$. Let $Z_t$ be the process defined in
(\ref{Z_process}). The same arguments as in the proof of Theorem
\ref{thm01} show that
\begin{eqnarray*}
E^x[t<T_n \wedge\tau_D;h(|X_t|^2)\varphi(\Phi_t)]&=&
 E^x[t<T_n \wedge\tau_D;h(|X_t|^2)\,
Z_t\,\exp(-\int_0^t q(R_s)ds)]\\
&=& C_k^{(\rho)}(1)E^{|x|^2}[t<T_n \wedge\tau_D;h(R_t)e_{-q}(t)]  \/.
\end{eqnarray*}
for every bounded Borel function $h$. We have, as before, $T_n \to \infty$
so letting $n \to \infty$ we obtain
\begin{equation*}
E^x[h(|X_t^D|^2)\varphi(\Phi_t^D))]=
C_k^{(\rho)}(1)E^{|x|^2}[h(R_t^{\widetilde{D}})e_{-q}(t)].
\end{equation*}
It means that
\begin{equation*}
(P_t^D\tilde{h}\tilde{\varphi})(x)=C_k^{(\rho)}(1)T_t^{\tilde{D}}h(|x|^2) \/,
\end{equation*}
where $\tilde{h}(x)=h(|x|^2)$,
$\tilde{\varphi}(y) =\tilde{\varphi}_x(y)=\varphi({(x,y) \o |x| |y|})$, and
 $\{T_t^{\tilde{D}}\}$ is the killed at $\partial \tilde{D}=\{r^2\}$
  Feynman-Kac semigroup based on the process $R_t$
with the potential $(-q)$.
Consequently,
\begin{equation*}
\int_0^\infty (P_t^D\tilde{h}\tilde{\varphi})(x)dt=
C_k^{(\rho)}(1)\,\int_0^\infty
T_t^{\tilde{D}}h(|x|^2)dt \/,
\end{equation*}
and
\begin{equation} \label{Green}
(G_D\tilde{h}\tilde{\varphi})(x)=C_k^{(\rho)}(1) \,V_{\tilde{D}} h(|x|^2) \/,
\end{equation}
where $G_D$ is the Green operator for the process $X_t$ and the set
$D$
and $V_{\tilde{D}}$ is the Green
 operator for the semigroup $\{T_t\}$ and the set $\tilde{D}$.

Integrating in polar coordinates we have

\begin{eqnarray} \label{Gpol}
(G_D\tilde{h}\tilde{\varphi})(x) &=&
 \int_D C_k^{(\rho)}(\cos \theta) G_D(x,y) h(|y|^2) dy \nonumber  \\
&=& \int_0^r h(R^2)\{\int_{S_R}C_k^{(\rho)}(\cos \theta) G_D(x,y)
d\sigma_R(y)\} dR \/,
\end{eqnarray}

and
\begin{equation} \label{Vpol}
V_{\tilde{D}} h(|x|^2) = \int_0^{r^2} h(y) V_{\tilde{D}}(|x|^2,y)dy
= \int_0^r h(R^2)V_{\tilde{D}}(|x|^2,R^2) 2R dR\/.
\end{equation}

Now, comparing \pref{Gpol}, \pref{Vpol},
\pref{Green} and applying the standard continuity arguments, we obtain


\begin{equation}
\label{greenequat01}
\widehat{(G_D)}_k(|x|,R) =
\frac{2}{\sigma_{n-1}R^{n-2}}\cdot V(|x|^2,R^2) \/.
\end{equation}
Since the $q$-Green function  $V_{\tilde{D}}(y,R^2)$ is $q$-harmonic in $y$,
 for $0<y<R$ and for $R<y<r$
(see \cite{ChZ}), we obtain
 that $\phi(y)=V_{\tilde{D}}(y,R^2)$ is a
solution of the same Schr\"odinger equation as in the case of the
gauge function in the proof of Theorem \ref{thm01} and the same computations
as before lead to the hypergeometric equation (\ref{equation3}).
Thus, on each of the intervals $(0,R)$, $(R,r]$, we can write
\begin{equation*}
V_{\tilde{D}}(|x|^2,R^2) = |x|^k\left(c_1(k,R^2)\,F_k(|x|^2) +
\frac{c_2(k,R^2)}{|x|^{2k+2\rho}} \, G_k(|x|^2)\right).
\end{equation*}
We have $\lim_{|x|^2\to r^2}V_{\tilde{D}}(|x|^2,R^2) = 0$. Therefore, for
$R<|x|<r$
\begin{equation*}
V_{\tilde{D}}(|x|^2,R^2) = c(k,R^2)\cdot|x|^k\left(\frac{G_k(|x|^2)}
{|x|^{2k+2\rho}}- \frac{G_k(r^2)}{r^{2k+2\rho}}\cdot
\frac{F_k(|x|^2)}{F_k(r^2)}\right).
\end{equation*}
As a consequence of the symmetry property of
$\widehat{(G_D)}_k(|x|,R)$ it follows easily that for $|x|<R$ we have
\begin{eqnarray*}
\widehat{(G_D)}_k(|x|,R) &=& \widehat{(G_D)}_k(R,|x|) \\
&=& \frac{2}{\omega_{n-1}|x|^{n-2}}\cdot V_{\tilde{D}}(R^2,|x|^2) \\
&=& \frac{2\,c(k,|x|^2)}{\omega_{n-1}|x|^{n-2}}\cdot
R^k\left(\frac{G_k(R^2)} {R^{2k+2\rho}}-
\frac{G_k(r^2)}{r^{2k+2\rho}}\cdot
\frac{F_k(R^2)}{F_k(r^2)}\right).
\end{eqnarray*}
To find $c(k,|x|^2)$ we use relation between Green function and
Poisson kernel.  Due to Green's theorem on hyperbolic space $\D^n$
we can obtain the Poisson kernel by differentiating the Green
function in the normal direction (see \cite{C}, page 174, theorem
8), i.e. for any $y\in S_1$ we get
$\left.-\frac{d}{dR}\right\vert_{R=r} G_D(x,Ry)=
(1-r^2)^{n-2}P_r(x,ry)$. Note that the factor $(1-r^2)^{n-2}$
appears because we consider the Poisson kernel as a density with
respect to the Lebesque measure $\sigma_r$ (not with respect to the
hyperbolic measure on the sphere). By bounded convergence theorem
and (\ref{pkgegenbauer}) we obtain
\begin{eqnarray*}
C_k^{(\rho)}(1)\,\omega_{n-1}\cdot\left(\left.-\dfrac{d}{dR}\right\vert_{R=r}\right)
\widehat{(G_D)}_k(|x|,R)
&=&\left(\left.-\dfrac{d}{dR}\right\vert_{R=r}\right)\int_{S_1}
C_k^{(\rho)}(\cos\theta)G_D(x,Ry)d\sigma_1(y)
\\
&=&\int_{S_1}\left(\left.-\dfrac{d}{dR}\right\vert_{R=r}\right)
C_k^{(\rho)}(\cos\theta)G_D(x,Ry)d\sigma_1(y)
\\
&=&(1-r^2)^{n-2}\int_{S_1}C_k^{(\rho)}(\cos\theta)P_r(x,ry)d\sigma_1(y)
\\
&=& \frac{(1-r^2)^{n-2}}{r^{n-1}}\int_{S_r}
C_k^{(\rho)}(\cos\theta)P_r(x,y) d\sigma_r(y)
\\
&=& \frac{(1-r^2)^{n-2}}{r^{n-1}}\,\, C_k^{(\rho)}(1)\cdot
\frac{|x|^kF_k(|x|^2)}{r^kF_k(r^2)}.
\end{eqnarray*}
On the other hand, using (\ref{hypgeo01}), we have
\begin{equation*}
-\left.\dfrac{d}{dR}\right\vert_{R=r} \left(\frac{G_k(R^2)}
{R^{2k+2\rho}}- \frac{G_k(r^2)}{r^{2k+2\rho}}\cdot
\frac{F_k(R^2)}{F_k(r^2)}\right) =
\frac{2(1-r^2)^{n-2}}{r^{n-1}}\cdot\frac{k+\rho}{r^{2k}\,F_k(r^2)}
\/.
\end{equation*}
Thus
\begin{equation*}
\frac{2\,c(k,|x|^2)}{\omega_{n-1}|x|^{n-2}} = C_n
\cdot\frac{\rho}{k+\rho}\cdot|x|^{k}\,F_k(|x|^2).
\end{equation*}
The proof is completed.
\end{proof}

\begin{thm}[Green function formula]
For $|x|<|y|$ we have
\begin{equation}
\label{gfformula}
 G_D(x,y) = C_n\sum_{k=0}^\infty
|x|^kF_k(|x|^2)|y|^k\left(\dfrac{G_k(|y|^2)}{|y|^{2k+2\rho}}
-\dfrac{G_k(r^2)}{r^{2k+2\rho}}\cdot\frac{F_k(|y|^2)}{F_k(r^2)}\right)\,C_k^{(\rho)}(\cos\theta).
\end{equation}
\end{thm}
\begin{proof}
We will consider both functions in (\ref{gfformula}) as a function
of $y$ on the sphere $S_R$ where $R=|y|$.  We now prove, under the
hypothesis $|x|<|y|$, that the series is uniformly convergent on
$S_R$. Indeed, it is easy to see that there exist constants
$c_1=c_1(n,z)$ such that $|G_k(z)|\leq c_1$. Moreover, we already
know that $|F_k(z)|\leq c_2$ and $|F_k(R^2)/F_k(r^2)|\leq c_3$
where $c_2, c_3$ do not depend on $k$. Thus the norm of the series is
bounded by
\begin{equation*}
c\left[|y|^{2-n}\,\sum_{k=0}^\infty
\left(\frac{|x|}{|y|}\right)^k\,C_k^{(\rho)}(1)
+r^{2-n}\sum_{k=0}^\infty\ \left(\frac{|x||y|}{r^2}\right)^k\,
C_k^{(\rho)}(1) \right]< \infty \/.
\end{equation*}
Since the terms of the series are continuous functions of $y$, the
series represents continuous function on the sphere $S_R$ which is
AS. Its Gegenbauer coefficients are the same as the coefficients
of the Green function computed in Theorem \ref{thm03}. This
completes the proof.
\end{proof}
\section{Examples}
In this section we provide an explicit integral formula for the Poisson kernel,
under some additional conditions imposed on our hypergeometric functions $F_k$.
We collect these conditions as Conjecture \ref{conjecture01} and \ref{conjecture02}.
Checking the validity of these conjectures we compute the exact form of the
integral formula for the Poisson kernel on $\D^4$ and $\D^6$.
We also give an integral formula for the Green function in these two cases.

We will need the following representations for the classical Newtonian
potentials and Poisson kernels on $\R^n$:
\begin{align}
    \label{transformD2}
    \log|x-y|^{-1} &=\log|y|^{-1}+\frac{1}{2}\sum_{k=1}^{\infty}\left(\frac{|x|}{|y|}\right)^k\,
   C_k^{(0)}(\cos\theta) \/,& |x|<|y| \/, \\
   \label{transformDn}
   |x-y|^{2-n}&=
  |y|^{2-n}\sum_{k=0}^{\infty}\left(\frac{|x|}{|y|}\right)^k\,
   C_k^{(\rho)}(\cos\theta) \/,& |x|<|y|,\,\,& n\geq3 \/, \\
   \frac{r^2-|x|^2}{|x-y|^n}
   \label{transformpk}
   &=r^{2-n}\sum_{k=0}^\infty \frac{k+\rho}{\rho}\cdot
    \left(\frac{|x|}{r}\right)^k C_k^{(\rho)}(\cos\theta) \/, &
    |x|<r,\,\,&|y|=r \/.
\end{align}
It is easy to obtain the first two formulas using the generating
function for Chebyshev and Gegenbauer polynomials respectively.
The last one is the consequence of (\ref{transformDn}) and the
relation $\frac{k+\rho}{\rho}\,\,C_k^{(\rho)}(x)=
C_k^{(\rho+1)}(x)-C_{k-2}^{(\rho+1)}(x)$, where $C_k^{(\rho)}(x)$
denotes zero if $k$ is negative. When $n=2$ so $\rho=0$ we obtain
 $k\,C_k^{(0)}(x)= C_k^{(1)}(x)-C_{k-2}^{(1)}(x)$, with the same
  convention as before for the meaning of
 negative indices $k$. From these relations and
(\ref{transformD2}) we obtain
\begin{equation}
\label{transformD2_2}
4\log|x-y|^{-1}=4\log|y|^{-1}-\frac{|x|^2}{|y|^2}
-2\sum_{k=1}^\infty
\left(\frac{|x|}{|y|}\right)^k\left(\frac{|x|^2}{|y|^2}\frac{1}{k+2}-
\frac{1}{k}\right)C_k^{(1)}(\cos\theta) \/.
\end{equation}

To simplify our notation we introduce the following notation:
\begin{equation*}
    f_z(k) = \frac{1}{\Gamma(k+{n\o2})}F_k(z) \,.
\end{equation*}
It is well-known that this is an entire function of complex variable
$k$. Observe that
\begin{equation*}
    \frac{F_k(|x|^2)}{F_k(r^2)} = \frac{f_{|x|^2}(k)}{f_{r^2}(k)}
\end{equation*}
for $k=0,1,2,\ldots$. We also define the function $H_{|x|^2,r^2}(k)$
by the following equality
\begin{equation}
\label{conjH}
 \frac{k+\rho}{\rho}\,\frac{f_{|x|^2}(k)}{f_{r^2}(k)}
= \left(\frac{1-|x|^2}{1-r^2}\right)^{\rho}\left(\frac{k+\rho}
{\rho}-\frac{n}{2}\,\frac{r^2-|x|^2}{(1-r^2)(1-|x|^2)}\right)+
\frac{r^2-|x|^2}{(1-r^2)^{n-2}}H_{|x|^2,r^2}(k).
\end{equation}
To proceed further we require two essential properties pertaining this function.
First of all, the function  $H_{|x|^2,r^2}$ is supposed to have the following
property:
\begin{equation*}
H_{|x|^2,r^2}(k) = O(k^{-1}) \quad \text{when} \quad |k|\to\infty\,.
\end{equation*}
This property is a
consequence of the following asymptotic expansion for the function
$F_k$, which we state in the sequel as the first conjecture

\begin{conj}
\label{conjecture01}
\begin{equation}
\label{conj01}
 F(k,-\rho;k+{n\o2};z) =
(1-z)^{\rho}+\frac{{n\o2}\,\rho}{k+{n\o2}}\,z(1-z)^{\rho-1}+O(k^{-2})\,.
\end{equation}
\end{conj}

Observe that Conjecture \ref{conjecture01} holds true for even $n$ as a consequence
of the well-known relation for
hypergeometric functions
\begin{equation*}
F(\alpha,\beta;\gamma;z) =
(1-z)^{-\beta}F\left(\gamma-\alpha,\beta;\gamma;\frac{z}{z-1}\right)\,.
\end{equation*}
 The above relation for odd $n$ is satisfied only for $0\leq z<{1\o2}$
 and we do not know if the above conjecture is satisfied by the function
 $F_k$ for odd $n$.

 The next conjecture, together with the first one, enable us to invert, in the
 sense of Laplace transform, our function  $H_{|x|^2,r^2}$. We formulate it as follows:

\begin{conj}
\label{conjecture02}
$F_k(z)$ as a function of (complex) variable $k$ has
no zeros in the region $\{\Re (k)\geq-n/2-\varepsilon\}$,
with
$\varepsilon=\varepsilon(z,n)>0$.
 \end{conj}

We now assume in the sequel that these two conjectures hold true.
 Then we define a function $w_{|x|^2,r^2}(v)$ as an inverse
Laplace transform (see, e.g. \cite{F})  of the function $H_{|x|^2,r^2}(k)$ and
consequently we get
\begin{equation}
\label{conj03}
 H_{|x|^2,r^2}(k) = \int_0^\infty e^{-kv}w_{|x|^2,r^2}(v)dv
\end{equation}
for every complex $k$ such that $\Re (k)\geq-n/2-\varepsilon$. Now
we can prove the following integral formula.
\begin{thm}
\label{intrepr}
Under assumptions as in Conjecture \ref{conjecture01} and \ref{conjecture02}
we obtain

\begin{equation}
\label{intformula}
 P_r(x,y) =
\frac{\Gamma({n\o2})}{2\pi^{n\o2}r(1-r^2)^{n-2}}\,\frac{r^2-|x|^2}{|x-y|^n}\int_0^\infty
\frac{w_{|x|^2,r^2}(v)L(x,y,v)}{|xe^{-v}-y|^{n-2}}\,\,dv
\end{equation}
where
\begin{eqnarray*}
L(x,y,v) &=& |xe^{-{v/2}}-ye^{v/2}|^{2\rho}\cdot\rho\,h(x,y,v)-
|x-y|^2\left[|xe^{-{v/2}}-ye^{v/2}|^{2\rho}-|x-y|^{2\rho}\right]\\
h(x,y,v) &=& |xe^{-{v/2}}-ye^{v/2}|^2-|x-y|^2 =
|x|^2(e^{-v}-1)+r^2(e^v-1).
\end{eqnarray*}
\end{thm}
\begin{proof}
Using (\ref{pkformula}) and (\ref{conjH}) we get
\begin{equation*}
P_r(x,y) = \textbf{A}_{\bf{1}}-\textbf{A}_{\bf{2}} +\textbf{B}\/,
\end{equation*}
where

\begin{eqnarray*}
\textbf{A}_{\bf{1}}&=& \frac{\Gamma({n \o 2})}{2\pi^{n \o
2}r^{n-1}}\left(\frac{1-|x|^2}{1-r^2}\right)^\rho
\sum_{k=0}^\infty\frac{k+\rho} {\rho}\frac{|x|^k}{r^k}
C_k^{(\rho)}(\cos\theta)\\
\textbf{A}_{\bf{1}}&=&\frac{\Gamma({n \o 2})}{2\pi^{n \o
2}r^{n-1}}\,\frac{n(r^2-|x|^2)(1-|x|^2)^{\rho-1}}{2(1-r^2)^{\rho+1}}\sum_{k=0}^\infty
\frac{|x|^k}{r^k} C_k^{(\rho)}(\cos\theta) \\
\textbf{B}&=&\frac{\Gamma({n \o 2})}{2\pi^{n \o
2}r^{n-1}}\,\frac{r^2-|x|^2}{(1-r^2)^{n-2}}\sum_{k=0}^\infty
\frac{|x|^k}{r^k}H_{|x|^2,r^2}(k)C_k^{(\rho)}(\cos\theta)\/.
\end{eqnarray*}
To deal with the series $\textbf{A}_{\bf{1}}$ and
$\textbf{A}_{\bf{2}}$ we use (\ref{transformpk}) and
(\ref{transformDn}) respectively and obtain
\begin{eqnarray*}
    \textbf{A}_{\bf{1}} &=&
    \frac{\Gamma({n \o 2})}{2\pi^{n \o
    2}r} \left(\frac{1-|x|^2}{1-r^2}\right)^\rho\cdot\frac{r^2-|x|^2}{|x-y|^n} \/, \\
    \textbf{A}_{\bf{2}} &=&
    \frac{\Gamma({n \o 2})}{2\pi^{n \o
    2}r}\cdot\frac{n}{2}\,\frac{(r^2-|x|^2)(1-|x|^2)^{\rho-1}}{(1-r^2)^{\rho+1}}\cdot\frac{1}{|x-y|^{n-2}} \/.
\end{eqnarray*}
We write $\textbf{B}$, using (\ref{conj03}), the Fubini's theorem
and (\ref{transformDn}), as
\begin{eqnarray*}
\textbf{B} &=& \frac{\Gamma({n \o 2})}{2\pi^{n \o
2}r^{n-1}}\cdot\frac{r^2-|x|^2}{(1-r^2)^{n-2}}\sum_{k=0}^\infty\frac{|x|^k}{r^k}
\left(\int_0^\infty
e^{-kv}w_{|x|^2,r^2}(v)dv\right)\,C_k^{(\rho)}(\cos\theta)\\
&=& \frac{\Gamma({n \o 2})}{2\pi^{n \o
2}r^{n-1}}\cdot\frac{r^2-|x|^2}{(1-r^2)^{n-2}}\int_0^\infty
\sum_{k=0}^\infty\frac{|xe^{-v}|^k}{r^k} \,C_k^{(\rho)}(\cos\theta)\,w_{|x|^2,r^2}(v)dv\\
&=&\frac{\Gamma({n \o 2})}{2\pi^{n \o
2}r}\cdot\frac{r^2-|x|^2}{(1-r^2)^{n-2}}\int_0^\infty
\frac{w_{|x|^2,r^2}(v)dv}{|xe^{-v}-y|^{n-2}}.
\end{eqnarray*}
Summing up $\textbf{A}_{\bf{1}}-\textbf{A}_{\bf{2}}+\textbf{B}$ we
obtain
\begin{eqnarray*}
P_r(x,y) &=&
\frac{\Gamma({n\o2})(r^2-|x|^2)}{2\pi^{n\o2}r(1-r^2)^{n-2}}
\left[\frac{(1-|x|^2)^\rho(1-|y|^2)^\rho}{|x-y|^n}-
\frac{n(1-|x|^2)^{\rho-1}(1-|y|^2)^{\rho-1}}{2|x-y|^{n-2}}\right.\\
&&\left.+ \int_0^\infty\frac{e^{\rho v} w_{|x|^2,r^2}(v)dv}
{|xe^{-{v\o2}}-ye^{v\o2}|^{n-2}}\right]\/.
\end{eqnarray*}
Recall that $h(x,y,v) = |x|^2(e^{-v}-1)+r^2(e^v-1)$. Using
(\ref{conj03}) and (\ref{conjH}) we get
\begin{eqnarray*}
    \int_0^\infty e^{\rho v} w_{|x|^2,r^2}(v)dv &=&
    H(-\rho) =
    \frac{n}{2}\,(1-|x|^2)^{\rho-1}(1-r^2)^{\rho-1}\\
    \int_0^\infty e^{\rho v}h(x,y,v)w_{|x|^2,r^2}(v)dv &=& |x|^2
    H(-\rho+1)+r^2H(-\rho-1)+(|x|^2+r^2)H(-\rho)\\
    &=& \frac{(1-|x|^2)^\rho(1-r^2)^\rho}{\rho}\/,
\end{eqnarray*}
where $H=H_{|x|^2,r^2}$ is the function defined in (\ref{conjH}). In
the last equality we use the following fact
\begin{equation*}
|x|^2\frac{f_{|x|^2}(-\rho+1)}{f_{r^2}(-\rho+1)} =
r^2\frac{f_{|x|^2}(-\rho-1)}{f_{r^2}(-\rho-1)}
\end{equation*}
which is a consequence of
\begin{equation*}
zf_z(-\rho+1) = z\sum_{i=0}^\infty
\frac{(-\rho+1)_i(-\rho)_i}{\Gamma(i+1)\Gamma(i+2)}z^i =
\frac{1}{\rho(\rho+1)}\sum_{i=0}^\infty
\frac{(-\rho)_{i+1}(-\rho-1)_{i+1}}{\Gamma(i+1)\Gamma(i+2)}z^{i+1} =
\frac{f_{z}(-\rho-1)}{\rho(\rho+1)}\/.
\end{equation*}
Thus finally we get
\begin{eqnarray*}
P_r(x,y) &=&
\frac{\Gamma({n\o2})}{2\pi^{n\o2}}\,\frac{r^2-|x|^2}{r(1-r^2)^{n-2}}
\left[\int_0^\infty \frac{\rho e^{\rho v}
h(x,y,v)w_{|x|^2,r^2}(v)dv}{|x-y|^n} - \int_0^\infty \frac{e^{\rho
v}w_{|x|^2,r^2}(v)dv}{|x-y|^{n-2}}\right.\\
&&\left.+ \int_0^\infty\frac{e^{\rho v} w_{|x|^2,r^2}(v)dv}
{|xe^{-{v\o2}}-ye^{v\o2}|^{n-2}}\right]\\
&=&
\frac{\Gamma({n\o2})}{2\pi^{n\o2}r(1-r^2)^{n-2}}\,\frac{r^2-|x|^2}{|x-y|^n}
\int_0^\infty \frac{e^{\rho v}
L(x,y,v)w_{|x|^2,r^2}(v)dv}{|xe^{-{v\o2}}-ye^{v\o2}|^{n-2}}\\
&=&\frac{\Gamma({n\o2})}{2\pi^{n\o2}r(1-r^2)^{n-2}}\,\frac{r^2-|x|^2}{|x-y|^n}
\int_0^\infty \frac{L(x,y,v)w_{|x|^2,r^2}(v)dv} {|xe^{-v}-y|^{n-2}}.
\end{eqnarray*}
This completes the proof.
\end{proof}
\begin{cor}
If $n=4$ then $w_{|x|^2,r^2}(v)=\frac{2(1+r^2)}{1-r^2}\, e^{-{2v\o
1-r^2}}$ and
\begin{equation*}
    P_r(x,y) =
    \frac{1}{2\pi^2r(1-r^2)^2}\cdot\frac{r^2-|x|^2}{|x-y|^4}\int_0^\infty
    \frac{[|x|^2(e^{-v}-1)+r^2(e^v-1)]^2\,w_{|x|^2,r^2}(v)}{|xe^{-v}-y|^2}\,dv.
\end{equation*}
If $n=6$ then $w_{|x|^2,r^2}(v)
=w^1_{|x|^2,r^2}(v)+w^2_{|x|^2,r^2}(v)$, where
\begin{eqnarray*}
    w^1_{|x|^2,r^2}(v) &=&
    3\left[1-r^2|x|^2+3(r^2-|x|^2)\right]\,e^{-bv}\cosh(cv)\\
    w^2_{|x|^2,r^2}(v) &=& -3\left[1+|x|^2r^2+5(|x|^2+r^2)\right]\frac{e^{-bv}\sinh(cv)}{2c}
\end{eqnarray*}
and
\begin{equation*}
    L(x,y,v) =
    [|x|^2(e^{-v}-1)+r^2(e^v-1)]^2(2|xe^{-{v\o2}}-ye^{v\o2}|^2+|x-y|^2),
\end{equation*}
hence
\begin{equation*}
    P_r(x,y) = \frac{1}{\pi^3r(1-r^2)^4}\cdot
    \frac{r^2-|x|^2}{|x-y|^6}\int_0^\infty
    \frac{L(x,y,v)w_{|x|^2,r^2}(v)\,dv}{|xe^{-v}-y|^4},
\end{equation*}
where $b=\frac{7-r^2}{2(1-r^2)}$,
$c=c(r)=\frac{\sqrt{r^4-14r^2+1}}{2(1-r^2)}$.
\end{cor}
\begin{proof}
For $n=4$ we have $F_k(z)=1-\frac{k}{k+2}\,z = (1-z)+\frac{2z}{k+2}$
(see Preliminaries) and it is obvious that Conjecture \ref{conjecture01} is valid
in this case. We also have
\begin{equation*}
    (k+1)\frac{F_k(|x|^2)}{F_k(r^2)} = (k+1)\,\frac{2+k(1-|x|^2)}{2+k(1-r^2)}
\end{equation*}
and it is analytic function of variable $k$ in the region $\{\Re
(k)\geq-2-\varepsilon\}$ for some $\varepsilon>0$ because the
denominator has only one zero $k_0 = -{2\o1-r^2}$. Thus we can use
Theorem \ref{intrepr}. We have
\begin{eqnarray*}
    L(x,y,v) = [|x|^2(e^{-v}-1)+r^2(e^v-1)]^2
\end{eqnarray*}
and all we have to do is to find the function $w_{|x|^2,r^2}(v)$. We
compute it using the residue method. We have
\begin{eqnarray*}
    \text{Res}_{k_0}H&=& \lim_{k\to k_0} (k-k_0)(k+1)\frac{F_k(|x|^2)}{F_k(r^2)}\\
    &=& \lim_{k\to k_0}(k-k_0)(k+1)\,\frac{2+k(1-|x|^2)}{2+k(1-r^2)}\\
    &=& \frac{r^2-|x|^2}{1-r^2}\cdot 2\frac{1+r^2}{1-r^2}
\end{eqnarray*}
and consequently
\begin{equation*}
    w_{|x|^2,r^2}(v) = 2\,\frac{1+r^2}{1-r^2}\,e^{-{2v\o1-r^2}}
\end{equation*}

For $n=6$ we have $F_k(z) = 1-\frac{2k}{k+3}\,z+ \frac{k(k+1)}
{(k+3)(k+4)}\,z^2$ and consequently
\begin{equation*}
\frac{k+2}{2}\,\frac{F_k(|x|^2)}{F_k(r^2)} =
\frac{k+2}{2}\cdot\frac{k(k+1)(1-|x|^2)^2+6k(1-|x|^2)+12}{k(k+1)(1-r^2)^2+6k(1-r^2)+12}
\/.
\end{equation*}
Zeros of the dominator are given by $k_1=-b-c$ and $k_2=-b+c$. It is
easy to check that Conjecture \ref{conjecture02} is valid also in
this case. We have
\begin{equation*}
    L(x,y,v)  =
    [|x|^2(e^{-v}-1)+r^2(e^v-1)]^2(2|xe^{-{v\o2}}-ye^{v\o2}|^2+|x-y|^2).
\end{equation*}
Observe that
\begin{eqnarray*}
    \text{Res}\,_{k_1}H+\text{Res}\,_{k_2}H
    &=& \frac{(1-r^2)^4}{r^2-|x|^2}\left[\lim_{k\to k_1}(k-k_1)\frac{k+2}{2}\,\frac{F_k(|x|^2)}{F_k(r^2)}+\lim_{ k\to
    k_2}(k-k_2)\frac{k+2}{2}\,\frac{F_k(|x|^2)}{F_k(r^2)}\right]\\
    &=& 3(1-|x|^2r^2+3(r^2-|x|^2))\/,\\
    \text{Res}\,_{k_1}H-\text{Res}\,_{k_2}H
    &=& \frac{(1-r^2)^4}{r^2-|x|^2}\left[\lim_{k\to k_1}(k-k_1)\frac{k+2}{2}\,\frac{F_k(|x|^2)}{F_k(r^2)}-\lim_{ k\to
    k_2}(k-k_2)\frac{k+2}{2}\,\frac{F_k(|x|^2)}{F_k(r^2)}\right]\\
    &=&-\frac{3}{c}\left[1+|x|^2r^2+5(|x|^2+r^2)\right].
\end{eqnarray*}
Thus we get
\begin{eqnarray*}
    w_{|x|^2,r^2}(v) &=& \frac{1}{2}\left(\text{Res}\,_{k_1}H+\text{Res}\,_{k_2}H\right)
    e^{-bv}\cosh(cv)
    +\frac{1}{2}\left(\text{Res}\,_{k_1}H-\text{Res}\,_{k_2}H\right)
    e^{-bv}\sinh(cv)\\
    &=& w^1_{|x|^2,r^2}(v)+w^2_{|x|^2,r^2}(v).
\end{eqnarray*}
The proof is completed.

Note that the polynomial $r^4-14r^2+1$, which appears in the
definition of the function $c(r)$, has zero in the interval $(0,1)$.
More precisely we have
\begin{align*}
    r^4-14r^2+1&>0, &0<r^2<7-4\sqrt{3} \/, \\
    r^4-14r^2+1&=0, &r^2=7-4\sqrt{3}  \/,\\
    r^4-14r^2+1&<0, &7-4\sqrt{3}<r^2<1 \/.
\end{align*}
If $r^2=7-4\sqrt{3}$ we have $c(r)=0$ and consequently
\begin{equation*}
w_{|x|^2,r^2}(v) = \left[1-r^2|x|^2+3(r^2-|x|^2)\right]e^{-bv}
-\left[1+|x|^2r^2+5(|x|^2+r^2)\right]\frac{ve^{-bv}}{2} \/.
\end{equation*}
In the case $7-4\sqrt{3}<r^2<1$ the function $c(r)$ has pure
imaginary values, i.e. $c(r)=i\tilde{c}(r)$ where
\begin{equation*}
    \tilde{c}(r) =\frac{\sqrt{|r^4-14r^2+1|}}{2(1-r^2)}  \/.
\end{equation*}
Thus we have
\begin{equation*}
w_{|x|^2,r^2}(v) =
\left[1-r^2|x|^2+3(r^2-|x|^2)\right]e^{-bv}\cos(\tilde{c}v)
-\left[1+|x|^2r^2+5(|x|^2+r^2)\right]\frac{e^{-bv}\sin(\tilde{c}v)}{2\tilde{c}}
\/.
\end{equation*}
\end{proof}
In the following corollary we provide the integral formula for the
Green function in $\D^4$.
\begin{cor}\label{cor01}
    For n=4 we get
    \begin{eqnarray}
        \label{greenD4}
        G_D(x,y) &=&
        \dfrac{1}{4\pi^2}\left(
        (1-|x|^2)(1-|y|^2) \left[\frac{1}{|x-y|^2}-\frac{r^2}{|y|^2|x-y^{*}|^2}\right]
        \right.\\
        \nonumber
        &&-4\left[\log\frac{1}{|x-y|}-\log\frac{r}{|y||x-y^{*}|}\right]
        + \left.(r^2-|x|^2)(r^2-|y|^2)
        \int_0^\infty\frac{w_{|x|^2,r^2}(v)\,dv}{|y|^2|xe^{-v}-y^{*}|^2}\right) \/,
    \end{eqnarray}
    where $w_{|x|^2,r^2}(v) = 2\frac{1+r^2}{1-r^2}\,e^{-{2v\o1-r^2}}$ and $y^*=\frac{r^2y}{|y|^2}$.
\end{cor}
\begin{proof}
For $n=4$ we have $G_0(z)=\frac{1}{z}+z+2\log z$ and
$G_k(z)=1-\frac{k+2}{2}\,z$. Thus, using (\ref{gfformula}) for
$n=4$, we obtain
\begin{equation*}
\frac{G_D(0,y)}{C_4}=\frac{G_0(|y|^2)}{|y|^2}-\frac{G_0(r^2)}{r^2}=\frac{1}{|y|^2}-|y|^2+4\log|y|
- \frac{1}{r^2}+r^2-4\log r \/.
\end{equation*}
where $C_4 = \frac{1}{4\pi^2}$. Elementary algebraic computation
shows that
\begin{eqnarray}
\label{greenD4_01} \frac{G_D(0,y)}{C_4} &=&
\left(\frac{(1-|x|^2)(1-|y|^2)}{|y|^2}+\frac{|x|^2}{|y|^2}+4\log|y|\right)
\\
\nonumber &&-
\left(\frac{(1-|x|^2)(1-|y|^2)}{r^2}+\frac{|x|^2|y|^2}{r^4}+4\log
r-\frac{1+r^2}{r^4}(r^2-|x|^2)(r^2-|y|^2)\right)\\
\nonumber &=& \bf{c}-\bf{d} \/.
\end{eqnarray}
 For $n=4$ the formula
(\ref{gfformula}) takes now the form
\begin{eqnarray*}
\frac{G_D(x,y)}{C_4} &=&
\frac{G_D(0,y)}{C_4}+\frac{1}{|y|^2}\sum_{k=1}^\infty
\left(\frac{|x|}{|y|}\right)^k(1-\frac{k}{k+2}|x|^2)(1-{{k+2}\o
k}|y|^2)C_k^{(1)}(\cos\theta)\\
&& -\frac{1}{r^2}\sum_{k=1}^\infty\left(\frac{|x||y|}{r^2}
\right)^k (1-\frac{k}{k+2}|x|^2)(1-{{k+2}\o k}r^2)\frac{(1-{k \o
{k+2}}|y|^2)} {(1-{k \o {k+2}}r^2)} C_k^{(1)}(\cos\theta) \\
&=& \frac{G_D(x,y)}{C_4} +\textbf{C} - \textbf{D} =
 (\bf{c}+\textbf{C})-(\bf{d}+\textbf{D}) \/.
\end{eqnarray*}
Put
\begin{eqnarray*}
   f_{|y|^2,r^2}(k) &=& \frac{1-{k \o {k+2}}|y|^2} {1-{k \o
    {k+2}}r^2} \/,\\
    g_{|x|^2,|y|^2}(k) &=& (1-\frac{k}{k+2}|x|^2)(1-{{k+2}\o k}|y|^2) \/.\\
\end{eqnarray*}
Then we have the following relations
\begin{eqnarray*}
    f_{|y|^2,r^2}(k) &=&
    \frac{1-|y|^2}{1-r^2}-2\,\frac{r^2-|y|^2}{(1-r^2)^2}\cdot\frac{1}{k+{2\o1-r^2}} \/;\\
    g_{|x|^2,|y|^2}(k) &=&
    (1-|x|^2)(1-|y|^2) +2\frac{|x|^2}{k+2}-2\frac{|y|^2}{k} \/;\\
    g_{|x|^2,r^2}(k)f_{|y|^2,r^2}(k) &=&
    (1-\frac{k}{k+2}|x|^2)(1-{{k+2}\o k}r^2)\cdot\frac{1-{k \o
    {k+2}}|y|^2}{1-{k \o {k+2}}r^2}\/.
\end{eqnarray*}
Observe also that the right-hand side of the last equality is equal
to
\begin{eqnarray*}
    &&(1-|x|^2)(1-|y|^2)+\frac{2|x|^2f_{|y|^2,r^2}(-2)}{(k+2)}-\frac{2r^2f_{|y|^2,r^2}(0)}{k}
    -2\,\frac{(r^2-|y|^2)}{(1-r^2)^2}\frac{g_{|x|^2,r^2}({-2\o1-r^2})}{k+{2\o1-r^2}}\\
    &&= (1-|x|^2)(1-|y|^2)
    +2\frac{|x|^2|y|^2}{r^2(k+2)}-2\frac{r^2}{k}
    - \frac{2(1+r^2)}{r^2(1-r^2)}\,
    \frac{(r^2-|x|^2)(r^2-|y|^2)}{k+{2\o{1-r^2}}}\/.
\end{eqnarray*}
Thus we get
$(\bf{c}+\textbf{C})=\textbf{C}_{\textbf{1}}+\textbf{C}_{\textbf{2}}$
and
$(\bf{d}+\textbf{D})=\textbf{D}_{\textbf{1}}+\textbf{D}_{\textbf{2}}-\textbf{D}_{\textbf{3}}$
where
\begin{eqnarray*}
\textbf{C}_{\textbf{1}} &=&
\frac{(1-|x|^2)(1-|y|^2)}{|y|^2}\sum_{k=0}^\infty\left(\frac{|x|}{|y|}\right)^k
C_k^{(1)}(\cos\theta) = \frac{(1-|x|^2)(1-|y|^2)}{|x-y|^2} \/,
\\
 \textbf{C}_{\textbf{2}} &=& 4\log|y|+\frac{|x|^2}{|y|^2}
+2\sum_{k=1}^\infty
\left(\frac{|x|}{|y|}\right)^k\left(\frac{|x|^2}{|y|^2}\frac{1}{k+2}-
\frac{1}{k}\right)C_k^{(1)}(\cos\theta) = -4\log\frac{1}{|x-y|}.
\end{eqnarray*}
The first equality follows from (\ref{transformDn}) and the other
is just (\ref{transformD2_2}). Similarly
\begin{eqnarray*}
\textbf{D}_{\textbf{1}} &=&
\frac{(1-|x|^2)(1-|y|^2)}{r^2}\sum_{k=0}^\infty
\left(\frac{|x||y|}{r^2}\right)^k C_k^{(1)}(\cos\theta)  =
\frac{(1-|x|^2)(1-|y|^2)r^2}{|y|^2|x-y^*|^2}  \/,
\\
\textbf{D}_{\textbf{2}} &=& 4\log
r+\frac{|x|^2|y|^2}{r^4}+2\sum_{k=1}^\infty
\frac{|x|^k|y|^k}{r^{2k}}\left(\frac{|x|^2|y|^2}{r^{4}}\frac{1}{k+2}-
\frac{1}{k}\right)C_k^{(1)}(\cos\theta) = -4\log\frac{r}{|y||x-y^*|}
\/.
\end{eqnarray*}
To deal with $\textbf{D}_{\textbf{3}}$ recall that $\int_0^\infty
e^{-{2v\o1-r^2}}e^{-kv}dv = \frac{1}{k+{2\o1-r^2}}$. Consequently
\begin{eqnarray*}
\textbf{D}_{\textbf{3}} &=&
\frac{2(1+r^2)}{r^4(1-r^2)}(r^2-|x|^2)(r^2-|y|^2)
\sum_{k=0}^\infty\left(\frac{|x||y|}{r^2} \right)^k
\frac{1}{k+{2\o{1-r^2}}}\,C_k^{(1)}(\cos\theta)\\
&=&
\frac{2(1+r^2)}{r^4(1-r^2)}(r^2-|x|^2)(r^2-|y|^2)\sum_{k=0}^\infty
\left(\frac{|x||y|}{r^2} \right)^k \bigl(\int_0^\infty
e^{-{2v\o1-r^2}}e^{-kv}dv \bigr)\,
C_k^{(1)}(\cos\theta)\\
&=&\frac{2(1+r^2)}{r^4(1-r^2)}(r^2-|x|^2)(r^2-|y|^2)\int_0^\infty\sum_{k=0}^\infty\left(\frac{|xe^{-v}||y|}{r^2}
\right)^k \,C_k^{(1)}(\cos\theta)e^{-{2v\o1-r^2}}dv\\
&=& 2\frac{1+r^2}{1-r^2}\,(r^2-|x|^2)(r^2-|y|^2)\int_0^\infty
\frac{e^{-{2v\o1-r^2}}dv}{|y|^2|xe^{-v}-y^*|^2} \/.
\end{eqnarray*}
Now summing up all the components we finally obtain (\ref{greenD4}).
\end{proof}
It is also possible to obtain an integral formula for the Green
function for $n=6$.
\begin{cor}
    For $n=6$ we have
    \begin{eqnarray*}
        G_D(x,y) &=&\frac{1}{4\pi^3}\left(
        (1-|x|^2)^2(1-|y|^2)^2\left[\frac{1}{|x-y|^4}-\frac{r^4}{|y|^4|x-y^*|^4}\right]\right.\\
        &&-6(1-|x|^2)(1-|y|^2)\left[\frac{1}{|x-y|^2}-\frac{r^2}{|y|^2|x-y^*|^2}\right]\\
        &&+24\left[\log\frac{1}{|x-y|}-\log\frac{r}{|y||x-y^*|}\right]-12\frac{(r^2-|x|^2)(r^2-|y|^2)}{|y|^2|x-y^*|^2}\\
        &&\left.+6(r^2-|x|^2)(r^2-|y|^2)\int_0^\infty
        \frac{W_{|x|^2,|y|^2}(v)dv}{|y|^4|xe^{-v}-y^*|^4}\right)\/,
    \end{eqnarray*}
    where
    \begin{eqnarray*}
         W_{|x|^2,|y|^2}(v) &=& f_1(x,y)e^{-bv}\cosh(cv)+f_2(x,y)\frac{e^{-bv}\sinh(cv)}{2(1-r^2)c}
    \end{eqnarray*}
    and
    \begin{eqnarray*}
        f_1(x,y) &=& r^2\frac{1+r^2}{1-r^2}(1-|x|^2)(1-|y|^2)+2(1-|x|^2|y|^2)-2(1-r^4)\\
        f_2(x,y) &=&
        r^2\frac{(1+r^2)^2}{1-r^2}(1-|x|^2)(1-|y|^2)+2(1-5r^2-2r^4)(1-|x|^2|y|^2)-2(1-r^2)^3\/.
    \end{eqnarray*}
\end{cor}
\begin{proof}
    We will show only some main steps of the proof of the above
    formula and explaining all details are left to the reader.
    Recall that for $n=6$ we have
    \begin{equation*}
        F_k(z) = 1-\frac{2k}{k+3}\,z+ \frac{k(k+1)}{(k+3)(k+4)}\,z
    \end{equation*}
    and
    \begin{equation*}
    G_k(z) = 1-\frac{2(k+4)}{k+1}\,z+\frac{(k+3)(k+4)}{k(k+1)}\,z\/.
    \end{equation*}
    Thus we get
    \begin{eqnarray*}
        F_k(|x|^2)G_k(|y|^2) &=&
        (1-|x|^2)^2(1-|y|^2)^2-6(1+|x|^2)(1+|y|^2)\left(\frac{|y|^2}{k+1}-\frac{|x|^2}{k+3}\right)\\
        &&+12\left(\frac{|y|^4}{k}-\frac{|x|^4}{k+4}\right)\\
        &=&(1-|x|^2)^2(1-|y|^2)^2-6(1-|x|^2)(1-|y|^2)\left(\frac{|y|^2}{k+1}-\frac{|x|^2}{k+3}\right)\\
        &&+12\left(\frac{|y|^4}{k(k+1)}-\frac{2|x|^2|y|^2}{(k+1)(k+3)}+\frac{|x|^4}{(k+3)(k+4)}\right)
    \end{eqnarray*}
    and
    \begin{eqnarray*}
    F_k(|x|^2)G_k(r^2)\frac{F_k(|y|^2)}{F_k(r^2)} &=&
    (1-|x|^2)^2(1-|y|^2)^2\\
    &-&6(1-|x|^2)(1-|y|^2)\left[\frac{r^2}{k+1}-\frac{|x|^2|y|^2}{r^2}\frac{1}{k+3}\right]\\
    &+&12\left(\frac{r^4}{k(k+1)}-\frac{2|x|^2|y|^2}{(k+1)(k+3)}+\frac{|x|^4|y|^4}{r^4}\frac{1}{(k+3)(k+4)}\right)\\
    &+&12(r^2-|x|^2)(r^2-|y|^2)\left[\frac{1}{k+1}-\frac{|x|^2|y|^2}{r^4(k+3)}\right]\\
    &+&\frac{6}{r^4}(r^2-|x|^2)(r^2-|y|^2)
    J(k)\/,
    \end{eqnarray*}
    where
    \begin{eqnarray*}
        J(k) &=&
        \left[f_1(x,y)\frac{1}{2}\left(\frac{1}{k+b+c}+\frac{1}{k+b-c}\right)
        +\frac{f_2(x,y)}{2c(1-r^2)}\frac{1}{2}\left(\frac{1}{k+b+c}-\frac{1}{k+b-c}\right)\right]\/.
    \end{eqnarray*}
    Observe also that
    \begin{eqnarray*}
        \frac{1}{|x-y|^2}&=& \frac{1}{|y|^2}\sum_{k=1}^\infty\left(\frac{|x|}{|y|}\right)^k
        \left(\frac{1}{k+1}-\frac{|x|^2}{|y|^2}\frac{1}{k+3}\right)C_k^{(2)}(\cos\theta)
        \\
        4\log|x-y|^{-1}&=&4\log|y|^{-1}-\frac{2}{3}\frac{|x|^2}{|y|^2}-\frac{|x|^4}{6|y|^4}\\
        &+&
        2\sum_{k=1}^\infty\frac{|x|^k}{|y|^k}\left[\frac{1}{k(k+1)}-\frac{|x|^2}{|y|^2}\frac{2}{(k+1)(k+3)}+
        \frac{|x|^4}{|y|^4}\frac{1}{(k+3)(k+4)}\right]C_k^{(2)}(\cos\theta)
    \end{eqnarray*}
    and
    \begin{eqnarray*}
        \frac{1}{r^8}\sum_{k=0}^\infty\frac{|x|^k|y|^k}{r^{2k}}\frac{1}{2}\left(\frac{1}{k+b+c}+\frac{1}{k+b-c}\right)
        C_k^{(2)}(\cos\theta)
        &=&\int_0^\infty
        \frac{e^{-bv}\cosh{(cv)}dv}{|y|^4|xe^{-v}-y^*|^4}\\
        \frac{1}{r^8}\sum_{k=0}^\infty\frac{|x|^k|y|^k}{r^{2k}}\frac{1}{2}\left(\frac{1}{k+b+c}-\frac{1}{k+b-c}\right)
        C_k^{(2)}(\cos\theta)
        &=&\int_0^\infty
        \frac{e^{-bv}\sinh{(cv)}dv}{|y|^4|xe^{-v}-y^*|^4}\/.
    \end{eqnarray*}
    Therefore using (\ref{gfformula}) for $n=6$ we can obtain the
    desired formula.
\end{proof}

\end{document}